\documentclass[12pt]{amsart}
\usepackage{amssymb}
\usepackage[dvips]{graphics}
\ifx\mathrm\undefined\let\mathrm\rm\fi
\ifx\mathbf\undefined\let\mathbf\bf\fi
\ifx\mathfrak\undefined\let\mathfrak\frak\fi
\ifx\mathcal\undefined\let\mathcal\cal\fi
\ifx\mathbb\undefined\let\mathbb\Bbb\fi
\ifx\emph\undefined\let\emph\it\fi
\font\bb=msbm10 at9.98pt
\begin{document}
\def\semidirect{\hbox{$\;$\bb\char'156$\;$}}
\newcommand{\SL}{\mathrm{SL}}
\newcommand{\GL}{\mathrm{GL}}
\newcommand{\g}{{{\mathfrak g}\,}}
\newcommand{\bor}{{{\mathfrak b}}}
\newcommand{\n}{{{\mathfrak n}}}
\newcommand{\h}{{{\mathfrak h\,}}}
\newcommand{\Id}{{\operatorname{Id}}}
\newcommand{\Z}{{\mathbb Z}}
\newcommand{\ZZ}{{\mathbb Z_{>0}}}
\newcommand{\N}{{\mathbb N}}
\newcommand{\R}{{\mathbb R}}
\newcommand{\p}{{\mathbb P}}
\newcommand{\C}{{\mathbb C}}
\newcommand{\Q}{{\mathbb Q}}
\newcommand{\CC}{\mathcal{C}}
\newcommand{\A}{\mathcal{A}}
\newcommand{\F}{\mathcal{F}}
\newcommand{\W}{\mathcal{W}}
\newcommand{\PP}{\mathcal{P}}
\newcommand{\Sym}{{\rm Sym}}
\newcommand{\Sing}{{\rm Sing}}
\newcommand{\Poly}{{\rm Poly}}
\newcommand{\Span}{{\rm Span}}
\newcommand{\Res}{{\rm Res}} 
\newcommand{\1}{{\bf 1}}
\newcommand{\kk}{{\bf k}}
\newcommand{\z}{{\bf z}}
\newcommand{\dontprint}[1]
{\relax}
\newtheorem%
{thm}{Theorem}
\newtheorem%
{prop}
{Proposition}
\newtheorem%
{lemma}
{Lemma}
\newtheorem%
{lemmadef}[thm]{Lemma-Definition}
\newtheorem%
{cor}
{Corollary}
\newtheorem%
{conj}
{Conjecture}
\newenvironment{definition}
{\noindent{\bf Definition\/}:}{\par\medskip}

\title {Asymptotic solutions to the $sl_2$ KZ equation\\
 and  the intersection of  Schubert classes}

\author[{}]
{I. Scherbak}

\maketitle     

\medskip
\centerline{\it School of Mathematical Sciences,
Tel Aviv University,}
\centerline{\it Ramat Aviv, Tel Aviv 69978, Israel}
\centerline{\it e-mail: \quad scherbak@post.tau.ac.il}

\bigskip

\pagestyle{myheadings}
\markboth{I. Scherbak}
{The $sl_2$ KZ equation and Schubert classes}
\begin{abstract}
${}$
The hypergeometric solutions to the KZ equation contain
a certain  symmetric ``master function'',  \cite{SV}. 
Asymptotics of the solutions correspond to critical points of
the master function and  give Bethe vectors of the  inhomogeneous Gaudin model, \cite{RV}. 
The general conjecture is that the  number of orbits of critical points equals the dimension 
of the relevant vector space,  and that the Bethe vectors form a basis. 
In \cite{ScV}, a proof of  the conjecture  for the $sl_2$ KZ equation was given.
The difficult part of the proof  was to count the number of  orbits
of critical points of the master function. 
 Here we present another, ``less technical'', proof  based on  a relation  
 between the master function and the map sending a set of polynomials
 into the  Wronski determinant. 
 Within these  frameworks, the number of  orbits  becomes the intersection number 
 of  appropriate Schubert classes. Application of the Schubert calculus  
to the $sl_p$  KZ  equation is discussed.

\end{abstract}
\section{Introduction}\label{S1}
Denote 
$Z=\left\{\ z=(z_1,\dots,z_n)\in\C^n\ \vert\  z_i\neq  z_j\,,\ 1\leq i<j\leq n\ \right\}\,$.
Let   $M=(m_1,\dots,m_n)$  be  a vector   with positive integer coordinates,
$|M|=m_1+\dots+m_n\,$. For  $z\in Z$, consider the following function in $k$ variables,  
$k\leq |M|/2\,$,
\begin{equation}
\label{Phi}
\Phi(t)=\Phi_{k,n}(t;z,M)=\prod_{i=1}^k\prod_{l=1}^n (t_i-z_l)^{-m_l}
\prod_{1\leq i<j\leq k}(t_i-t_j)^2\,, 
\end{equation}
defined on
\begin{eqnarray}
\label{CC}
&\ &\CC(t,z)\,=\, \CC_{k,n}(t,z)\,=\\
&= &\, \left\{\, t=(t_1,\dots,t_k)\in\C^k\ \vert\
t_i\neq z_j,\ t_i\neq t_l,\ 1\leq i<l\leq k,\ 1\leq j\leq n\, \right\}\,.\nonumber
\end{eqnarray}
The critical point system of the function   $\Phi(t)\,$
provides the Bethe equations for the $sl_2$ Gaudin model of an inhomogeneous magnetic 
chain,  \cite{G}.
The symmetric group acts on $\CC(t,z)$ permuting coordinates $t_1,\dots,t_k$,
and the action preserves the critical set of  $\Phi(t)$. 
The orbits of nondegenerate critical points 
give common  eigenvectors of the  system
of commuting Hamiltonians of the Gaudin model and define asymptotic solutions
to the $sl_2$ KZ equation.
For generic $z\in Z$, these eigenvectors  
generate the  subspace  $\Sing_k$  of singular vectors of  weight  $|M|-2k$ 
in the tensor product of irreducible $sl_2(\C)$ representations 
with highest weights $m_1,\dots,m_n\,$, \cite{RV}.

\medskip
One of the main results in  \cite{ScV} asserts that for generic $z$ 
all critical points of the function   
$\Phi(t)$ are nondegenerate and the number of orbits  equals the dimension of   $\Sing_k$.  
In other words, each common eigenvector  of the Hamiltonians of the Gaudin model  is 
represented by the solutions to the Bethe equations exactly once. 
In general, for  other integrable 
models  ``parasite''  solutions to the Bethe equations appear.

\medskip
The difficult part of the proof given in \cite{ScV} is to get an 
appropriate upper bound for the number of  orbits 
of critical points of the function $\Phi(t)$. The aim of the present paper  is 
to calculate this upper bound in another, ``less technical'', 
way using the Schubert calculus.

\medskip
The Wronski determinant of  two polynomials  in one variable defines 
a map from the Grassmannian of two-dimensional  planes of the linear
space of polynomials  to the space of monic polynomials (Sec.~\ref{s31}).
It turns out that  there is a one-to-one correspondence between the
orbits of critical points of  the  function $\Phi(t)$  and the planes
in the preimage under this map of the polynomial
\begin{equation}
\label{W}
W(x)=(x-z_1)^{m_1}\dots (x-z_n)^{m_n}\,.
\end{equation}
In fact, this is a classical result going back to Heine and Stieltjes
(Sec.~\ref{s33}; cf.\cite{S}, Ch.~6.8). 
  
\medskip
To calculate the cardinality of the preimage is a problem of enumerative 
algebraic geometry, and an upper bound can be easily obtained in terms of
the intersection  number of special Schubert classes (Sec.~\ref{S4}). 
A well-known relation between  representation theory and  the Schubert calculus
implies that the obtained upper bound coincides with the dimension  
of the space $\Sing_k\,$.

\medskip
The Schubert calculus seems to be a useful tool for the $sl_p$ KZ equations
as well. A function related  in a very similar way  to the Wronski map from 
the Grassmannian of $p$-dimensional  planes was found by  A.~Gabrielov, 
\cite{Ga}.  This function turns to be the  master function of the
 $sl_p$ KZ equation associated  with the tensor product of  symmetric powers 
of the standard $sl_p$-module. The Schubert calculus provides an upper bound
for the number of orbits of critical points as the intersection number of appropriate
Schubert classes.  Arguments similar to those of the $sl_2$ case prove then
that the Bethe vectors form a basis of the corresponding space of singular vectors,  
Sec.~\ref{s51}. 

\medskip
The Wronski map corresponds  
to some special rational curve in the Grassmannian (\cite{EG}, see also
Sec.~\ref{s52}). It would be interesting to find  rational curves 
corresponding to the master functions associated with other tensor products.

\medskip
 Our  presentation  of the KZ equation  theory  in Sec.~\ref{S2} is borrowed from  
 \cite{RV}, \cite{V}.   
The  exposition  of  the Schubert calculus and  its relation to the  
$sl_2$ representation 
theory  in Sec.~\ref{S3} and \ref{S4} follows   \cite{GH} and \cite{F}.

\medskip
The author is thankful to A.~Eremenko,  A.~Gabrielov, F.~Sottile and A.~Varchenko 
for useful  suggestions, to G.~Fainshtein and D.~Karzovnik for stimulating comments.

\section{The master function}\label{S2}
\subsection{$sl_2$ modules}\label{s21}

The Lie algebra $sl_2=sl_2(\C)$ is generated by
 elements $e, f, h$ such that $[e,f]=h,\  [h,e]=2e,\ [h,f]=-2f\,.$
     
\medskip
For   a nonnegative integer  $m$, denote $L_{m}$  the irreducible 
$sl_2$-module with  highest weight $m$ and  highest weight singular 
vector   $v_m$. We have $hv_m=mv_m$,\ $ev_m=0$. 

 \medskip
 Denote $S_m$   the unique bilinear symmetric form  on $L_{m}$ such that
$S_m(v_m,v_m)=1$ and  $S_m(hx,y)=S_m(x,hy)\,,\  S_m(ex,y)=S_m(x,fy)$ for all
$x,y\in L_{m}$. The   vectors 
$v_m, fv_m,\dots, f^{m}v_m$ are orthogonal with respect to $S_m$ and
form a basis of  $L_{m}$.

\medskip
 The tensor product 
\begin{equation}
\label{LM}
L^{\otimes M}=L_{m_1}\otimes \dots \otimes  L_{m_n}\,,
\end{equation}
where  $M=(m_1,\dots, m_n)$ and all $m_i$ are nonnegative integers,  
becomes an $sl_2$-module  if the action of $x\in sl_2$ is defined by
$$
x \otimes 1\otimes \dots \otimes 1+1\otimes x\otimes \dots \otimes 1
\dots + 1\otimes 1\otimes \dots \otimes x.
$$
The  bilinear symmetric form  on  $L^{\otimes M}$  given by  
 \begin{equation}
 \label{S}
S=S_{m_1}\otimes\dots\otimes S_{m_n}
\end{equation}
is  called {\it the Shapovalov form}.

\medskip
Let $j_1,\dots,j_n$ be  integers such that $0\leq j_i\leq m_i$
for any $1\leq i\leq n$. Denote
\begin{equation}\label{fJ}
f^Jv_M=f^{j_1}v_{m_1}\otimes \dots \otimes f^{j_n}v_{m_n}\,,\ J=(j_1,\dots,j_n)\,.
\end{equation}
The vectors $\left\{f^Jv_M\right\}$ are orthogonal with respect to the Shapovalov  
form  and provide a basis of  the space $L^{\otimes M}$.

\medskip
For  $A=(a_1,\dots,a_n)$,  write $|A|=a_1+\dots +a_n$.  We have 
$$
h(f^Jv_M)\,=\,(|M|-2|J|)f^Jv_M\,,\quad e(f^Jv_M)=0\,,
$$
i.e. the vector $f^Jv_M$ is {\it a singular vector of weight} $|M|-2|J|$.

\subsection{The $sl_2$ KZ equation}\label{s22}

For the Casimir element,                                                    
$$
\Omega=e\otimes f+f\otimes e+\frac12 h\otimes h \in sl_2\otimes sl_2,
$$
and for $1\leq i<j\leq n$, denote $\Omega_{ij}:L^{\otimes M}\to L^{\otimes M}$  
the operator which acts  as $\Omega$ on $i$-th and $j$-th factors
and as the identity on all others. For  $z\in Z$, define operators 
$H_1(z),\dots, H_n(z)$ on $L^{\otimes M}$ as follows,
\begin{equation}
\label{H}
H_i(z)=\sum_{j\neq i}\frac{\Omega_{ij}}{z_i-z_j}\,, \quad  i=1,\dots, n.
\end{equation}

\medskip
{\it The  Knizhnik--Zamolodchikov (KZ) equation} on a function 
$u: Z\to L^{\otimes M}$  is the system of partial differential equations
\begin{equation}
\label{KZ}
\kappa \frac{\partial u}{\partial z_i}=H_i(z)u\,,\quad   i=1,\dots, n,
\end{equation}         
where $\kappa$ is a parameter. This equation appeared first 
in Wess--Zumino models of conformal field theory, \cite{KZ}.

\subsection{Asymptotic solutions to the KZ equation {\rm (\cite{V})}}\label{s23}
\medskip
The series
\begin{equation}
\label{phi}
\phi(z)=e^{I(z)/\kappa}\left(\phi_0(z)+\kappa \phi_1(z)+\kappa^2 \phi_2(z)+
\dots \right),
\end{equation}
where $I$ and $\phi_i$ are  functions on $Z$, provides
{\it an asymptotic solution to the equation}  (\ref{KZ}) if the substitution                             
of $\phi$ into (\ref{KZ}) gives $0$ in the
expansion into a formal power series in $\kappa$.   

\medskip
If $\phi(z)$ is  an asymptotic solution, then the  substitution  into   (\ref{KZ})  
gives  recurrence relations for the functions $\phi_i$, 
and in general  $\phi_i$ can be recovered  from $\phi_0$
for any $i>0$. Moreover the substitution shows that
\begin{equation}\label{ev}
 \phi_0(z)\ \  {\it is\  an\  eigenvector\  of\    }  H_i(z)\  \ {\it with\  the\  
eigenvalue\  } \frac{\partial I}{\partial z_i} (z)
\end{equation}
for any $1\leq i\leq n$.

\subsection {The Bethe Ansatz}\label{s24}

For any $z\in Z$, the operators  $H_i(z)$ commute, are symmetric with respect 
to the Shapovalov form,  and therefore have a common basis of eigenvectors.

\medskip
{\it The Bethe Ansatz}  is a method for diagonalizing of a system 
of commuting operators $H_i$ in a vector space $L$ called {\it the space of states}. 
The idea is to consider some  $L$-valued function $w\,:\,\C^k\,\to\, L$, 
and to determine its argument
$t=(t_1,\dots,t_k)$ in such a way  that the value of this function, 
$w(t)$, is an eigenvector.
The  equations which determine these special values of the argument
are called {\it the Bethe equations}
and the vector  $w(t^0)$, where $t^0$ is a solution to the Bethe equations,  
is called {\it the Bethe vector}.

\medskip
In all known examples   the  Bethe equations coincide with the critical point
system  of a suitable function  called  
{\it the master function} of the model.  The Bethe vectors correspond
to the critical points of   the master function. 

\medskip                           
The standard conjectures are that the (properly counted) number of critical 
points of the master function is equal to the dimension of the space of states, 
and that the Bethe vectors form a basis.

\medskip
The operators $H_i(z)$ given by (\ref{H}) appear in the  $sl_2$ Gaudin
model of inhomogeneous magnetic chain \cite{G}.
 According to (\ref{ev}), the asymptotic solutions to the  $sl_2$ KZ equation
 (\ref{KZ}) are labeled by the Bethe vectors of the  $sl_2$ Gaudin model.

\subsection{Subspaces of singular vectors}\label{s25}

Let $k$ be a nonnegative integer, $k\leq  |M|/2$, 
and $\Sing_k\,=\,\Sing_k (L^{\otimes M})$ the subspace of 
singular vectors of weight  $|M|-2k$ in $L^{\otimes M}$,
\begin{equation}
\label{Sing}
\Sing_k (L^{\otimes M})=\left\{\, w\in  L^{\otimes M}\ \vert\  ew=0\,,\ hw=(|M|-2k)w\,
\right\}.
\end{equation}
This space is generated by the vectors $f^Jv_M$  as in (\ref{fJ}) with $|J|=k$.

\begin{thm}  {\rm (\cite{ScV}, Theorem~5 )}\label{tI} We have
$$
\dim\Sing_k(L^{\otimes M})\,=\,\sum_{q=0}^{n} (-1)^q
\sum_{1\leq i_1<\dots <i_q\leq n}
{k+n-2-m_{i_1}-\dots -m_{i_q}-q\choose n-2}\,.
$$
\end{thm}
As usually we set ${a\choose b}=0$ for $a<b$.

\medskip
The number $\dim\Sing_k(L^{\otimes M})$ is clearly the multiplicity of 
$L_{|M|-2k}$ in the decomposition of the tensor product
$L^{\otimes M}$ into a direct sum of irreducible  $sl_2$-modules.
Equivalently, this is the multiplicity of the trivial $sl_2$-module,
$L_0\,$, in the decomposition of the tensor product
$$
L_{m_1}\,\otimes\, \dots\, \otimes \, L_{m_n}\, \otimes \, L_{|M|-2k}\,.
$$
It seems the explicit formula was unknown to experts in representation theory.

\medskip
The KZ equation preserves  $\Sing_k $ for any $k$,  and the
subspaces $\Sing_k $    generate
the whole of $L^{\otimes M}$. Therefore it is enough to produce solutions
with values in a given $\Sing_k=\Sing_k (L^{\otimes M})$. 

\subsection{Hypergeometric solutions  {\rm (\cite{SV})}}\label{s26}  

For $z\in Z$,  consider the function
 $$
\Psi_{k,n}(t,z)=\Psi_{k,n}(t,z;M)=\prod_{1\leq i<j\leq n}(z_i-z_j)^{m_im_j/2}
\prod_{i=1}^k\prod_{l=1}^n (t_i-z_l)^{-m_l}
\prod_{1\leq i<j\leq k}(t_i-t_j)^2 \,,
$$
which is defined on  $\CC(t,z)$  given by  (\ref{CC}). 
This function is invariant with respect to the group $S^k$
of  permutations of $t_1,\dots,t_k$.

\medskip
For $J=(j_1,\dots,j_n)$ with integer coordinates satisfying 
$0\leq j_i\leq m_i$ and  $|J|=k$,  set
$$
A_J(t,z)=\frac1{j_1!\dots j_n!}\Sym_t\left[\, \prod_{l=1}^n
\prod_{i=1}^{j_l}\frac1{t_{j_1+\dots+j_{l-1}+i}-z_l}\, \right]\,,
$$
where  
$$
\Sym_tF(t)=\sum_{\sigma\in S^k}F\left(t_{\sigma(1)},\dots, t_{\sigma(k)}\right)
$$
is the sum  over all  permutations  of $t_1,\dots, t_k$.

\begin{thm}\label{t1}   {\rm (\cite{SV})}\ \  
Let $\kappa\neq 0$ be fixed. The function
$$
u^{\gamma(z)}=\sum_{J\,:\,|J|=k}\left(\int_{\gamma(z)}\Psi^{1/\kappa}_{k,n}(t,z)
A_J(t,z)dt_1\wedge\dots \wedge dt_k\right) \cdot f^Jv_M\,,
$$
where   $\gamma(z)$ is  an appropriate   $k$-cycle lying  in $\CC_{k,n}(z)$,
provides a  nontrivial solution to the KZ equation  {\rm (\ref{KZ})} and
takes values in $\Sing_k$. 
\end{thm}

\subsection{The master function of the Gaudin model  
{\rm (\cite{RV})}}\label{s27}
Asymptotic solutions can be produced by taking the limit of $u^{\gamma(z)}$
as $\kappa \to 0$. According to the steepest descend method,  
the leading terms, $\phi_0$, are defined by the critical points of  
$\Psi^{1/\kappa}_{k,n}(t,z)$ with respect to $t$. In order to study the critical points, 
one can clearly replace $\Psi^{1/\kappa}_{k,n}(t,z)$ with the function  $\Phi(t)$ given by 
(\ref{Phi}). 

\begin{thm}\label{t2}{\rm (\cite{RV})}\ \   (i)\ \
If $t^0$ is a nondegenerate critical point of  the function $\Phi(t)$, then the vector
$$
w(t^0,z)=\sum_{J\, :\, |J|=k} A_J(t^0,z)f^Jv_M\,
$$
is an eigenvector of the set of commuting operators $H_1(z),\dots, H_n(z)$.

\noindent
(ii)\ \
For generic $z\in Z$, the eigenvectors $w(t^0,z)$ generate the space $\Sing_k\,$. 
\end{thm}

Thus the function  $\Phi(t)$ is  the master function of the Gaudin model,
and  $w(t^0,z)$ are the Bethe vectors. The critical point system of the
function  $\Phi(t)$ is as follows,
\begin{equation}\label{cps}
\sum_{l=1}^n\frac{m_l}{t_i-z_l}\,=\,\sum_{j\neq i}\frac{2}{t_i-t_j}\,,
\quad  i=1,\dots, n\,.
\end{equation}

\medskip
The set of  critical points of the function  $\Phi(t)$ is invariant
with respect to the permutations of $t_1,\dots,t_k$. 
Critical points belonging to the same 
orbit clearly define the same vector.
One of the main results of  \cite{ScV} is as follows.

\begin{thm}\label{t3}  {\rm (\cite{ScV})}\ \   For generic $z\in Z$,

\noindent
(i) \ \  all critical points of the function  $\Phi(t)$ given by 
{\rm  (\ref{Phi})}  are nondegenerate;

\noindent
(ii)\ \  the number of orbits of critical points equals the dimension of  
$\Sing_k\,$.
\end{thm}

Thus for generic $z\in Z$, the Bethe vectors form a basis of  $\Sing_k\,$.

\medskip
The  proof of  the statement {\it (i)} of Theorem \ref{t3} 
given in \cite{ScV} is short, see  Theorem~6  in  \cite{ScV}. 
The statements $(ii)$ of Theorem \ref{t2} and $(i)$ of Theorem 
\ref{t3} imply that the number of orbits is {\it at least} $\dim\Sing_k\,$.
The difficult part of  
\cite{ScV} was to prove that  the number of orbits is  
{\it at most} $\dim\Sing_k\,$, see Theorems 9--11 in  \cite{ScV}.

\medskip
As it was pointed out in \cite{ScV},  Sec.~1.4, the orbits of critical points 
of  the function  $\Phi(t)$ are labeled  by  certain two-dimensional 
planes in the linear space of complex polynomials. This observation 
suggests to apply the Schubert calculus to the problem under 
consideration. 

\section{The Bethe vectors and  the Wronski map}\label{S3}
\subsection{The Wronski map}\label{s31}
Let  $\Poly_d$ be  the vector space of  complex polynomials of degree $\leq d$
 in one variable and let  $G_2(\Poly_d)$ be the Grassmannian of  two-dimensional 
 planes in  $\Poly_d$. The complex dimension of    $G_2(\Poly_d)$ is $2d-2$.

 \medskip
 For any $V\in G_2(\Poly_d)$ define  {\it the degree of} 
of $V$ as  the maximal degree  of its polynomials and the {\it order}  of $V$ as  
the minimal   degree of its non-zero polynomials.
Let $V\in G_2(\Poly_d)$ be a plane of order $a$ and of degree $b$. 
Clearly $a<b$. Choose in $V$ two  monic polynomials, $F(x)$ and $G(x)$,
 of  degrees   $a$ and $b$ respectively.  They form  a basis of $V$.   
{\it The Wronskian of} $V$  is defined as the monic polynomial
$$
W_V(x)=\frac{F'(x)G(x)-F(x)G'(x)}{a-b}\,.
$$

The following lemma is evident.
\begin{lemma}\label{l1}
(i)\ \ The degree of $W_V(x)$ is  $a+b-1\leq 2d-2$.

\noindent
 (ii)\ \ The   polynomial  $W_V(x)$  does not depend on the  choice of a monic basis.

  \noindent
 (iii)\ \ All polynomials  of   degree $a$  in $V$ are proportional.
\hfill $\triangleleft$
\end{lemma}

 Thus  the mapping sending $V$ to $W_V(x)$ is  a well-defined map from   $G_2(\Poly_d)$
 to $\C\p^{2d-2}$. We call it  {\it the Wronski map}. This is a mapping between 
 smooth complex algebraic varieties of the same dimension,  and  hence the preimage
 of  any polynomial consists of  a finite number of   planes. 
On Wronski maps see \cite{EG}.
 
 \subsection{Planes with a given Wronskian}\label{s32}
\begin{lemma}\label{l2}
 Any plane with a given Wronskian is uniquely determined by any 
of its  polynomial. 
\end{lemma}

\noindent
{\bf Proof:}\  \  Let $W(x)$ be the Wronskian of a plane $V\,$ and  $f(x)\in V$.
Take any  polynomial  $g(x)\in V$ linearly independent with $f(x)$. 
The plane $V$  is the solution space of the following second order   linear
differential equation  with respect to  unknown function   $u(x)$,
$$
\left|\begin{array}{ccc}u(x)\ \ &f(x)\   \ &g(x)\\
                          u'(x)\ \ &f'(x)\ \  &g'(x)\\
                          u''(x)\  \ &f''(x)\ \ &g''(x)
			   \end{array}\right|=0.
$$
The Wronskian of the polynomials $f(x)$ and $g(x)$ is clearly proportional 
to $W(x)$. Therefore this equation can be re-written in the form 
$$
W(x)u''(x)-W'(x)u'(x)+h(x)u(x)=0\,,
$$
where 
$$
h(x)\, =\, \frac{-W(x)f''(x)+W'(x)f'(x)}{f(x)}\,,
$$ 
as $f(x)$ is clearly  a solution to this equation.
\hfill $\triangleleft$

\medskip
We  call a plane $V$  {\it generic} if for any $x_0\in\C$  there is a polynomial 
$P(x)\in V$ such that  $P(x_0)\neq 0$. In a generic plane, the polynomials of 
any basis do not have common roots, and  almost all polynomials  of the bigger degree 
do not have multiple roots. 
 
\medskip
The following lemma is evident.
\begin{lemma}\label{l3}  Let $V$  be  a generic plane.

 \noindent
(i)\ \  If  $P(x)=(x-x_1)\dots(x-x_l)\in V$ is a polynomial  without
multiple roots, then $W_V(x_i)\neq 0$ for all $1\leq i\leq l$. 

 \noindent
(ii)\ \  If   $x_0$  is   a root of multiplicity  $\mu>1$ of  a
polynomial  $Q(x)\in V$, then  $x_0$ is a root of  $W_V(x)$
of multiplicity  $\mu-1$.  \hfill $\triangleleft$
\end{lemma}

A generic plane $V$ is {\it nondegenerate} if  the polynomials  of  the smaller degree 
 in $V$ do not have multiple roots.

\begin{lemma}\label{l4} Let $V$ be a nondegenerate plane. 
If  the Wronskian $W_V(x)$ has the form 
$$
W_V(x)=x^{m}\tilde W(x)\,,\ \  \tilde W(0)\neq 0\,,
$$
then there exists a polynomial   $F_0(x)\in V$ of the form  
$$
F_0(x)=x^{m+1}\tilde F(x)\,,\ \  \tilde F(0)\neq 0\,.
$$
\end{lemma}

\noindent
{\bf Proof:}\ \  Let $G(x)\in V$ be a polynomial of the smaller degree.
We have $G(0)\neq 0$, due  to Lemma \ref{l3}. Let
 $F(x)\in V$  be a polynomial of the bigger degree. The polynomial
$$
F_0(x)\,=\,F(x)-\frac{F(0)}{G(0)}\,G(x)\,\in V
$$
satisfies  $F_0(0)=0$  and therefore  has the form
$$
F_0(x)=x^l  \tilde F(x)\,, \ \ \tilde F(0)\neq 0\,,
$$
for some integer  $l\geq 1$.  The polynomials $G$ and $F_0$ form a basis
of $V$, therefore  the polynomial $F'_0(x)G(x)-F_0(x)G'(x)\,$ is 
proportional to $W_V(x)$.  The smallest degree term in
this polynomial is  $l\cdot a_l\cdot G(0)x^{l-1}\,$,
where  $a_l$ is the coefficient of $x^l$ in $F_0(x)$.  Therefore $l=m+1$.
 \hfill $\triangleleft$

\subsection{ The master function and  the Wronski map}\label{s33}
In the XIX century, Heine and Stieltjes  in their studies of  second
order linear differential equations with polynomial coefficients and a polynomial
solution of a prescribed degree (cf. the proof of Lemma \ref{l2} in Sec.~\ref{s32})
arrived at the  result, which can be formulated as follows.

\begin{thm}\label{tW} {\rm (Cf. \cite{S},  Ch.~6.8)}\quad Let 
$t^0$  be a critical point of  the function $\Phi(t)$   given by {\rm (\ref{Phi})}.
Then $F(x)=(x-t^0_1)\dots(x-t^0_k)$ is  a polynomial of the smaller degree
in a nondegenerate plane with  the Wronskian $W(x)$ given by  {\rm (\ref{W})}.

Conversely,  if $F(x)=(x-t^0_1)\dots(x-t^0_k)$ is a  polynomial  of the smaller
degree in  a plane $V$ such that $W_V(x)=W(x)$, then 
 $t^0=(t^0_1,\dots, t^0_k)$ is a critical point of  $\Phi(t)$.
\hfill $\triangleleft$   
\end{thm}

\begin{cor}\label{c2}
There is a one-to-one correspondence 
between the  orbits of critical points of  the  master function   
$\Phi(t)$ given by {\rm (\ref{Phi})}  
and the nondegenerate planes  of  order $k$  and of degree $|M|+1-k$  having   
Wronskian $W(x)$ given by {\rm (\ref{W})}.
\end{cor}

\noindent
{\bf Remark.}\ \ The function  $\Phi(t)$   given by (\ref{Phi}) can be 
re-written in the form
$$
 \Phi(t)\,=\,\Phi (W\,, F)=\,\frac{{\rm Disc}^2(F)}{\Res(W,F)}\,,
$$
where $\Res(W,F)$ is the resultant of the polynomials $W(x)$ and
$F(x)$,  ${\rm Disc}(F)$ is the discriminant of the polynomial $F(x)$,
the polynomial $W(x)$ is given by (\ref{W})  and 
$F(x)=(x-t_1)\dots(x-t_k)$. 
If $W=F'Q-FQ'$ for some polynomial $Q=Q(x)$,
then $W'=F''Q-FQ''\,$, and we have 
$$
\frac{W'}{W}\,=\,\frac{F''Q-FQ''}{F'Q-FQ'}\,.
$$
If polynomials $F$ and $Q$ do not have common roots, then at each root $t_i$ 
of $F$ we clearly have
$$
\frac{W'(t_i)}{W(t_i)}\,=\,\frac{F''(t_i)}{F'(t_i)}\,,\quad  i=1,\dots,k.
$$
This is exactly the critical point system of the function $ \Phi(t)$, cf. (\ref{cps}).
Thus the function $ \Phi$  considered as a function of polynomials $W$ and $F$ is 
{\it the generating function} of the Wronski map -- for a given polynomial $W$, 
the critical points of the generating function determine the nondegenerate planes
in the preimage of $W$. 
 
\section{The preimage of a given Wronskian}\label{S4}

The number of  nondegenerate planes  with a given  Wronskian 
can be estimated from above by the intersection number of Schubert classes.

\subsection{Schubert calculus {\rm (\cite{GH}, \cite{F})}}\label{s41}
Let   $G_2(d+1)= G_2(\C^{d+1})$ be the Grassmannian variety of  
two-dimensional subspaces $V\subset \C^{d+1}$.   
A chosen basis  $e_1,\dots,e_{d+1}\,$ of $ \C^{d+1}$
defines the flag of linear subspaces
$$
E_{\bullet}\,:\quad E_1\, \subset \, E_2\, \subset\, \dots\,  
\subset \, E_d \subset \, E_{d+1}= \C^{d+1}\,,
$$
where $E_i=\Span\{e_1,\dots,e_i\}\,$, $\dim E_i=i$.   
For any integers  $a_1$ and $a_2$
such that   $0\leq a_2\leq a_1\leq d-1\,$, the {\it Schubert variety} 
$\Omega_{a_1,a_2}(E_{\bullet})\subset  G_2({d+1})$
is defined as follows,
$$
\Omega_{a_1,a_2}= \Omega_{a_1,a_2}(E_{\bullet})=\left\{\, V\in  G_2({d+1})\, 
\vert\, \dim \left(V\cap E_{d-a_1}\right)\geq 1\,,
\   \dim \left(V\cap E_{d+1-a_2}\right)\geq 2\, \right\}\,.
$$
The variety $\Omega_{a_1,a_2}=\Omega_{a_1,a_2}(E_{\bullet})$
is an irreducible closed subvariety of $ G_2({d+1})$  of the complex
codimension $a_1+a_2$.

\medskip
The homology classes $[\Omega_{a_1,a_2}]$  of Schubert varieties
 $\Omega_{a_1,a_2}$ are independent of the choice of  flag, and form a basis
for the integral homology of  $G_2(d+1)$.
Define  $\sigma_{a_1,a_2}$ to be the cohomology class in 
 $H^{2(a_1+a_2)}(G_2({d+1}))$ whose cap product with the
fundamental class of   $G_2(d+1)$ is the homology class
 $[\Omega_{a_1,a_2}]$.
The classes  $\sigma_{a_1,a_2}$ are called {\it Schubert classes}.
They give a basis over $\Z$ for the cohomology ring of the Grassmannian.
The product or {\it intersection} of  any two Schubert classes   
$\sigma_{a_1, a_2}$ and  $\sigma_{b_1, b_2}$
has the form
$$
\sigma_{a_1, a_2}\,\cdot\,\sigma_{b_1, b_2}\,=\,
\sum_{c_1+c_2=a_1+a_2+b_1+b_2} C(a_1,a_2;b_1,b_2; c_1,c_2) 
\sigma_{c_1, c_2}\,,
$$
where $ C(a_1,a_2;b_1,b_2; c_1,c_2)$ are nonnegative integers
called {\it the Littlewood--Richardson coefficients}.

\medskip
If the sum of the codimensions of classes equals $\dim G_2({d+1})=2d-2$,
then their intersection is an integer (identifying the generator of the top
cohomology group $\sigma_{d-1, d-1}\in H^{4d-4}(G_2({d+1}))$ with $1\in\Z$) 
called {\it the intersection number}.

\medskip
When $(a_1, a_2)=(q,0)$, $0\leq q\leq d-1$, 
the  Schubert varieties $\Omega_{q, 0}$  
are called {\it special} and the corresponding cohomology classes 
$\sigma_q\,=\,\sigma_{q, 0}$
are called {\it special Schubert classes}.

\medskip
We will need the following  fact connecting the 
Schubert calculus and representation theory.  

\medskip
Denote $L_q$ the irreducible $sl_2$-module with highest weight $q$.

 \begin{prop}\label{Pr} {\rm (\cite{F})}
 Let $q_1,\dots, q_{n+1}$ be integers such that   $0\leq q_i\leq  d-1$  
for all $1\leq i\leq n+1$ and  $q_1+\dots +q_{n+1}=2d-2$. Then the
intersection number of the corresponding special Schubert classes,
$ \sigma_{q_1}\cdot {\rm\ ...\ } \cdot \sigma_{q_{n+1}}\,$,  coincides with  
the multiplicity of  the trivial $sl_2$-module  $L_0$ in
the tensor product  $L_{q_1}\otimes\dots\otimes L_{q_{n+1}}$.  \hfill $\triangleleft$
 \end{prop}

This proposition and Theorem \ref{tI} imply the explicit formula for the
intersection number of special Schubert classes which
seems to be unknown to experts.

\begin{cor}  Let $q_1,\dots, q_{n+1}$ be integers such that   $0\leq q_i\leq  d-1$
for all $1\leq i\leq n+1$ and  $q_1+\dots +q_{n+1}=2d-2$. Then we have
\begin{eqnarray}\label{int}
\sigma_{q_1}\cdot {\rm\ ...\ } \cdot \sigma_{q_{n+1}}=\sum_{l=1}^n (-1)^{n-l}
\sum_{1\leq i_1<\dots <i_l\leq n}
{q_{i_1}+\dots +q_{i_l}+l-d-1\choose n-2}\,.
\end{eqnarray}
\end{cor}

\subsection{The  planes with a given Wronskian  and  Schubert classes}\label{s42}

Applying the Schubert calculus to the Wronski map, 
we arrive at the following result.

\begin{thm}\label{tWr}
Let $m_1\,,\dots, m_n\,$ be positive integers such that   
$|M|=m_1+\dots+m_n\leq 2d-2$. 
Then for generic $z\in Z$  and for any integer $k$ 
such that $0<k<|M|+1-k\leq d$,
the preimage of  the polynomial  {\rm (\ref{W})}  
under the Wronski map consists of at most
 $$
\sigma_{m_1}\,\cdot  {\rm\ ...\ } \cdot\, \sigma_{m_n}\,\cdot \,\sigma_{|M|-2k}
$$
planes of order $k$ and of degree $<|M|+1-k$.
\end{thm}

\noindent{\bf Proof:}\ \  Any plane of  order $k$ with the Wronskian of degree $|M|$
 lies in  $G_2( \Poly_{|M|+1-k})$, as Lemma \ref{l1} shows.
Therefore it is enough to consider the Wronski map on   $G_2( \Poly_{|M|+1-k})$.

\medskip
For  any  $z_i$, define the  flag  $\F_{z_i}$ in $ \Poly_{|M|+1-k}\,$, 
$$
 \F_0(z_i) \subset \F_1(z_i)\subset \dots \F_{|M|+1-k}(z_i),=\, \Poly_{|M|+1-k}\,,
$$
where  $\F_j(z_i)$  consists of the polynomials 
$ P(x)\in\Poly_{|M|+1-k}\,$ of the form
$$
 P(x)= a_j(x-z_i)^{|M|+1-k-j}+\dots +a_0(x-z_i)^{|M|+1-k}\,.
 $$
We have  $\dim\F_{j}(z_i)=j+1$.
\medskip
Lemma \ref{l4} implies that the nondegenerate planes  with a Wronskian
having  at $z_i$ a root of multiplicity $m_i$ lie in the special
Schubert variety   
$$
\Omega_{m_i,  0}(\F_{z_i})\subset G_2(\Poly_{|M|+1-k})\,.
$$

\medskip
The maximal possible degree of  the Wronskian  $W_V(x)$ for 
$V\in\Poly_{|M|+1-k}$ is clearly $2|M|-2k$.
Denote $m_{\infty}\,=\,|M|-2k$, the difference between $2|M|-2k$ and $|M|$.
If $m_{\infty}$ is positive, it is
{\it the multiplicity of $W(x)$ at infinity}. 

\medskip
The nondegenerate planes  with a Wronskian
having  given multiplicity $m_{\infty}$ at infinity lie in   
the special Schubert variety   
$\Omega_{m_{\infty}, 0}(\F_\infty)\,$, where $\F_\infty$ is the flag
$$
\Poly_0\subset\Poly_1\subset\dots\subset \Poly_{|M|+1-k}\,.
$$

\medskip
We conclude that the nondegenerate planes which have the Wronskian
 (\ref{W})  lie in the intersection of special 
Schubert varieties
$$
\Omega_{m_1, 0}(\F_{z_1})\,\cap \Omega_{m_2, 0}(\F_{z_2})
\cap \dots \cap \Omega_{m_n, 0}(\F_{z_n})
\, \cap \Omega_{m_{\infty}, 0}(\F_{\infty})\,.
$$
\medskip
The dimension of  $G_2(\Poly_{|M|+1-k})$  is exactly
$m_1+\dots+m_n+m_{\infty}\,$, therefore  this intersection 
consists of a finite number planes, and
the intersection number of the special Schubert classes 
$$
\sigma_{m_1}\,\cdot {\rm\ ...\ }\cdot \sigma_{m_n}\,\cdot\, \sigma_{|M|-2k}\,
$$
provides an upper bound.   \hfill $\triangleleft$    

\medskip
The statement $(ii)$ of Theorem \ref{t2}, Corollary \ref{c2},  
Proposition \ref{Pr} and Theorem \ref{tWr} prove the
statement  {\it (ii)} of Theorem \ref{t3}.

\begin{cor}
For generic $z\in Z$, all planes in the preimage under the Wronski map
of the polynomial $W(x)$ given by {\rm (\ref{W})}  are nondegenerate.
\end{cor}

\section{Comments: $sl_p$ case}\label{S5}
\subsection{The Schubert calculus and the $sl_p$ KZ equations}\label{s51}
Any irreducible  $sl_2$-module is a symmetric power of the standard one.
It turns out that our arguments  work for the $sl_p$  KZ equations 
associated  with the tensor product of symmetric powers of the standard 
$sl_p$-module. These arguments show that 

\medskip\centerline
{\it the number of orbits of critical points
of the master function}

\centerline
{\it equals the dimension of the relevant subspace
of singular vectors} 

\medskip\noindent
and hence prove the conjecture on Bethe vectors for 
the corresponding Gaudin model. Here we describe briefly the  
steps following the scheme of  $sl_2$ case (Sec.~\ref{S2}--\ref{S4}).

\bigskip\noindent 
(I)\ {\it The master function.}\quad
For  the Lie algebra $sl_p=sl_p(\C)$, denote $E$ the standard $sl_p$-module
and $\alpha_1,\dots, \alpha_{p-1}$  the simple roots.
Let $M=(m_1,\dots,m_n)$ be  a fixed vector  with positive integer coordinates.
Set   $L_j=\Sym^{m_j}E\,$ and $L=L_1\otimes\dots \otimes L_n\,$.
Every $L_j$ is the irreducible  $sl_p$-module with highest weight
$\Lambda_j=(m_j,0,\dots,0)$.  Write $\Lambda=\Lambda_1+\dots
+\Lambda_n\,$. Let  $k_1,\dots,k_{p-1}$ be integers,
$k_1\geq \dots \geq k_{p-1}\geq 0$,   $\kk=(k_1,\dots,k_{p-1})$. 
Denote $\Sing_{\kk} L$ the subspace 
of singular  vectors  in $L$ of weight  $\Lambda-k_1\alpha_1-\dots-
-k_{p-1}\alpha_{p-1}$.

\medskip 
For the KZ equation associated with  the subspace $\Sing_{\kk} L$,
the master function $\tilde\Phi (t)=\tilde\Phi_{K,n}(t;z;M)$ 
is  a function  in  $K=k_1+\dots+k_{p-1}$ complex variables  $t$, 
$$
t=\left(\, t^{(1)}\,,\     \dots\,,  t^{(p-1)}\,\right)\,,\quad   
t^{(l)}= \left(\, t_1^{(l)}\,,\  \dots\,,\  t_{k_l}^{(l)}\right)\,,
\quad  l=1,\dots, p-1\,,
$$
defined on $\CC_{K,n}(t;z)$, see (\ref{CC}),
and given by
\begin{eqnarray}\label{P}
\tilde\Phi (t)=\tilde\Phi_{K,n}(t;z;M)& \,=\,&\prod_{l=1}^{p-1}\ \ 
\prod_{1\leq i<j\leq k_l}\ \ 
(t^{(l)}_i-t^{(l)}_j)^2 \\
&\times &\prod_{l=2}^{p-1}\ \prod_{i=1}^{k_{l-1}}\ \prod_{j=1}^{k_l}
(t^{(l-1)}_i-t_j^{(l)})^{-1} \nonumber\\
&\times &\prod_{j=1}^n\ \prod_{i=1}^{k_1}\, 
(t^{(1)}_i-z_j)^{-m_j}\,,  \nonumber
\end{eqnarray}
according to \cite{RV}.  
The set of critical points of the function (\ref{P}) 
is invariant with respect to the group
$S=S^{k_1}\times\dots\times S^{k_{p-1}}\,$, 
where  $S^{k_l}$ is the group of 
permutations of $t_1^{(l)},\dots,t_{k_l}^{(l)}\,$.
If $k_l=0$ for  $l_0\leq l\leq p-1$, then
corresponding terms in (\ref{P}) are missing.

\bigskip\noindent 
(II)\ {\it The Wronski map.}\quad The Wronskian of  $l$ polynomials 
$f_1(x),\dots, f_l(x)\,$ is  the determinant
$$
{\rm Wr}[f_1,\dots,f_l](x)=\det
 \left(\begin{array}{ccc}f_1(x)&\dots &f_l(x)\\
                          f_1^\prime(x)&\dots &f_l^\prime(x)\\
                               \dots &\dots &\dots\\   
                               f_1^{(l-1)}(x)&\dots & f_l^{(l-1)}(x)
                          \end{array}\right)\,.
$$
Denote  $G_p(\Poly_d)$ the Grassmannian of   $p$-dimensional planes in $\Poly_d$.
The Wronskian of  any  $V\in G_p(\Poly_d)$  is defined as a monic polynomial 
which is proportional to the Wronskian  of some (and hence, any) basis of  $V$. 
{\it The Wronski map} \ $\W: \, G_p(\Poly_d)\, \to\, \C\p^{p(d+1-p)}$\ sends
$V\in G_p(\Poly_d)$ into its Wronskian, \cite{EG}. 

\medskip
Any element of $G_p(\Poly_d)$ has a basis of polynomials of pairwise distinct
degrees, and these degrees are uniquely defined  by  the element.
For $V\in G_p(\Poly_d)$, let  $P_1(x),\dots, P_p(x)$ be a basis such that
$\deg P_l=d_l$ for all $1\leq l\leq p$ and   $d\geq d_1 >d_2 >\dots > d_p\geq 0.$
Denote ${\rm Wr}_l(x)={\rm Wr}[P_{l+1},\dots, P_p](x)$ for any $0\leq l\leq p-1$.
In particular,   ${\rm Wr}_0(x)$ is proportional to the Wronskian of  $V$, 
and ${\rm Wr}_{p-1}(x)=P_p(x)$.  We have $\deg {\rm Wr}_l=k_l$, where
\begin{equation}\label{k}
 k_l=d_p+\dots+d_{l+1}-\frac{(p-l)(p-l-1)}2\,,\ \ 0\leq l\leq p-1.
\end{equation}
If   polynomials ${\rm Wr}_l(x)$  and   ${\rm Wr}_{l+1}(x)$
do not have common roots for all
$0\leq l \leq p-2$, we call   $V=\Span \{P_1(x),\dots, P_p(x)\}\in G_p(\Poly_d)$
{\it a nondegenerate $p$-plane of type} $D=(d_1,\dots, d_p)$.

\medskip
The question of enumerative algebraic geometry  we are interested in
is as follows:

\medskip\centerline
{\it Given monic polynomial $W(x)$ and  integers $d_1> \dots > d_p\geq 0$,}

\centerline 
{\it find the number of nondegenerate $p$-planes}

\centerline
{\it of type $D=(d_1,\dots, d_p)$ with the Wronskian $W(x)$.}

\medskip
Such $p$-planes  belong to $G_p(\Poly_{d_1})$, therefore we write $d=d_1$. 

\bigskip\noindent 
(III) {\it Critical points of the master function and
nondegenerate $p$-planes.}\quad 
Fix  integers $d_1>d_2 >\dots> d_{p}\geq 0$, $d_1=d$. 
Let $k_1,\dots, k_{p-1}$ be as in (\ref{k}). 
One can re-write the function (\ref{P})  in the  form  
(cf. Remark in Sec. \ref{s33})
$$
\tilde\Phi (t)\,=\, {{\rm Disc}^2(W_1)\cdots{\rm Disc}^2(W_{p-1})\over
{\rm Res}(W,W_1){\rm Res}(W_1,W_2)\cdots{\rm Res}(W_{p-2},W_{p-1})}\,,
$$
where  $W_l=W_l(x;t^{(l)})=\prod_{j=1}^{k_l}(x-t^{(l)}_j)$  
is a   polynomial in $x$ of degree  $k_l$  for  $1\leq l\leq p-1$ 
and $W=W(x)=\prod_{j=1}^n(x-z_j)^{m_j}$ is the same as in  (\ref{W}).
In this form the function $\tilde\Phi (t)$ was found by A.~Gabrielov 
as the generating function of the Wronski map, \cite{Ga}.
If degrees of polynomials $W_l$ vanish for ${l_0}\leq l\leq {p-1}$, 
then the corresponding terms in $\tilde\Phi(t)$ are missing.
The following statement is a generalization of  
the result of Heine and Stieltjes (cf. Theorem~\ref{tW} in
Sec.~\ref{s33}).

\begin{thm}\label{GS} 
Let {\rm (\ref{W})}  be the Wronskian of  some $p$-dimensional plane 
in the space of polynomials.  
Assume that there exists a basis $\{P_1(x),\dots, P_p(x)\}$ 
of this plane such that  the Wronskian of the polynomials 
$P_{l+1}(x),\dots, P_p(x)$ is $W_l(x;s^{(l)})=\prod_{j=1}^{k_l}(x-s_j^{(l)})\,$  
for any $1\leq l\leq p-1$ and $s=(s^{(1)},\dots,s^{(p-1)})\in \CC_{K,n}(t;z)$.
Then the point $s$ is a critical point  of the master function {\rm (\ref{P})}.

Conversely,  any orbit of critical points  of the  master function 
{\rm (\ref{P})}  defines a unique $p$-dimensional plane with these properties.
\end{thm}

\begin{cor}\label{cp}
The number of orbits of critical points of the master function
{\rm (\ref{P})} coincide with the number of nondegenerate $p$-planes of type
$D=(d_1,\dots, d_p)$ in the preimage of {\rm (\ref{W})} under the Wronski map.    
Here $d_1=d,\ d_l=k_{l-1}-k_l+p-l\  (2\leq l\leq p-1),\ d_p=k_{p-1}.$
\end{cor}

\bigskip\noindent 
(IV) {\it An upper bound: the Schubert calculus.}\quad
The nondegenerate  $p$-planes of Corollary \ref{cp} lie in the intersection
of  Schubert varieties
$$
\Omega_{(m_1)}(\F_{z_1})\,\cap \Omega_{(m_2)}(\F_{z_2})
\cap \dots \cap \Omega_{(m_n)}(\F_{z_n})
\, \cap \Omega_{w}(\F_{\infty})\,,
$$
where  $\Omega_{(m_j)}=\Omega_{m_j,0,\dots,0}$ are special Schubert
varieties and   $ \Omega_{w}= \Omega_{w_1,\dots,w_p}$ is a Schubert
variety with
$w_l=d+l-p-d_{p+1-l}$,\  $1\leq l\leq p$ or,  in terms of $\kk$,
$$
w_1=d+1-p-k_{p-1}\,,\quad w_l=d+1-p-k_{p-l}+k_{p-l+1}  {\rm \ for \ } 
2\leq l\leq p-1, \quad  w_p=0.
$$
Thus the number of orbits of critical points of the master function (\ref{P}) 
is bounded from above by the intersection number of the corresponding
Schubert classes. This number is the dimension of $\Sing_{\kk}L$, \cite{F}.  

\bigskip\noindent 
(V) {\it A lower bound: Fuchsian equations.}\quad
In the $sl_2$ case this step was done in \cite{RV} (see Sec.~\ref{s27}). 
In the $sl_p$ case, results of N.~Reshetikhin and A.~Varchenko 
(see Theorem~10.4 of \cite{RV}) say that the bound from below by 
the same number, $\dim \Sing_{\kk}L$, holds for generic $\z$ 
if this is the case  for $n=2, z_1=0, z_2=1$.  
This is really the case, as we will see 
using the theory of Fuchsian differential equations (\cite{R}).

\medskip 
Denote  $\tilde\Phi^0(t)=\tilde\Phi_{K,2}(t;\{0,1\};\{m_1,m_2\})$.
This is the master function of the $sl_p$ KZ equation associated with
the subspace of singular vectors of the weight
$\Lambda_1+\Lambda_2-k_1\alpha_1-\dots - k_{p-1}\alpha_{p-1}$ 
in the tensor product of two symmetric powers of the standard $sl_p$-module.

\medskip
 Assume $m_1\geq m_2$.
The Pieri formula (\cite{F}) says that the dimension of this subspace
is $1$ if $0\leq k_1\leq m_2$ and $k_2=\dots=k_{p-1}=0$ and $0$ otherwise.
 
\medskip
According to Theorem~\ref{GS}, a critical point of the function
$\tilde\Phi^0(t)$ corresponds to the solution space of the linear
differential equation of order $p$ with respect to  unknown function   
$u(x)$,
\begin{equation}\label{eW}
\det \left(\begin{array}{cccc} u(x)\ \ &P_1(x)\   \ &\dots\   \ &P_p(x)\ \ \\
                          u'(x)\ \ &P_1'(x)\ \  &\dots\ \ &P_p'(x)\    \\
                           \dots\ \  &\dots\ \   &\dots\ \  &\dots\ \\ 
                          u^{(p)}(x)\   &P_1^{(p}(x)\ \ &\dots\ \ &P_p^{(p)}(x)\\
                           \end{array}\right)=0,
\end{equation}
where $\{P_1(x),\dots,P_p(x)\}$ is a basis of this $p$-plane. 
This is a Fuchsian differential equation with regular singular points
at $0$, $1$ and $\infty$. 
The exponents at $0$ and  $1$ are
$0,1,\dots,p-2, m_1+p-1$ and $0,1,\dots,p-2, m_2+p-1$, respectively.
The Wronskian of this equation is the Wronskian of  the $p$-plane,
$$
{\rm Wr}[P_1,\dots, P_p](x)=x^{m_1}(x-1)^{m_2}.
$$

\begin{prop}\label{p2}
The equation {\rm (\ref{eW})} has the form
$$
x(x-1)u^{(p)}(x)+(Ax+B)u^{(p-1)}(x)+Cu^{(p-2)}(x)=0.
$$
\end{prop}

Hence we can assume that polynomials  $P_3(x),\dots, P_{p-1}(x), P_p(x)$ 
are $x^{p-3}, \dots, x,\,1 $, respectively.
Then the Wronskians $W_l(x)={\rm Wr}[P_{l+1},\dots, P_p](x)=1$ for 
$2\leq l\leq p-1$, i.e. $k_2=\dots=k_{p-1}=0$.

\begin{cor}
The function $\tilde\Phi^0(t)$ may have critical points only if  
 $k_2=\dots=k_{p-1}=0$. 
\end{cor}
One way to calculate the number of critical points of   the function 
$\tilde\Phi^0(t)$ is to note that for   $k_2=\dots=k_{p-1}=0$ 
it becomes the function $\Phi_{k_1,2}(t;\{0,1\};\{m_1,m_2\})$, 
see (\ref{Phi}).
Another way is to observe that the equation of Proposition~\ref{p2} 
is the hypergeometric equation  with respect to the function $u^{(p-2)}(x)$, \cite{R}.
In any way we get that the number of orbits of critical points of the function 
 $\tilde\Phi^0(t)$ is $1$ when the corresponding highest weight enters 
the tensor product, and $0$ otherwise. 

\medskip
 Thus the upper and lower bounds for the number of orbits of critical
points of the function $\tilde\Phi(t)$ coincide, and this number is
the dimension of $\Sing_{\kk}L$.

\subsection{Master functions and rational functions in the Grassmannian}\label{s52}
As it is explained in \cite{EG}, any rational curve of degree $A$ 
in the Grassmannian $G_p(\Poly_d)$ defines a projection from  
$G_p(\Poly_d)$ to the projective space of complex polynomials of degree
at most $A$, considered up to proportionality.  The construction is as follows.

\medskip
Let $V\in  G_p(\Poly_d)$.  One can identify polynomials  of  
degree at most $d$ with vectors in  
$\C^{d+1}$ using coefficients as coordinates. Then $V$ can be given by a  
$p\times (d+1)$-matrix
of coordinates  of  $p$ linearly independent  polynomials of $V$.  
On the other hand,  one can also identify these polynomials 
with linear forms on  $\C^{d+1}$, and 
$p$ linearly independent forms define a subspace $\tilde V$ of dimension  
$d+1-p$ in $\C^{d+1}$. 
This subspace can be represented by a  $(d+1-p)\times (d+1)$-matrix. 
We denote this matrix $K_V$.

\medskip
Now let a rational curve be given by a $p\times (d+1)$-matrix $P(\xi)$ 
of polynomials in $\xi$, and let  $K_V$  
be a $(d+1-p)\times (d+1)$-matrix corresponding in the described way  to   
any $V\in G_p(\Poly_d)$.
Then  $P(\xi)$  defines a projection by the formula  
\begin{equation}\label{map}
V\mapsto  \det  \left(\begin{array}{c}P(\xi)\\ K_V  \end{array}\right)\,.
\end{equation}
As it was shown in \cite{EG}, the Wronski map corresponds 
to a rational curve  given by the matrix
 $$
 \left(\begin{array}{c}F(\xi)\\ F^\prime(\xi)\\ \dots\\   
F^{(p-1)}(\xi) \end{array}\right)\,,
$$
where $F(\xi)$ is a row of powers of $\xi$,
$$
F(\xi)\, = \, (\xi^d\quad \xi^{d-1}\quad\dots \quad \xi \quad 1)\,.
$$
Hence this rational curve  corresponds to the master function (\ref{P}).
It would be interesting to understand what rational curves correspond  
to the master functions for  $sl_p$ KZ equations associated with other 
tensor products.


\newpage
\bigskip

\begin{thebibliography}{1} 

\bibitem[EG]{EG} A.~Eremenko and A.~Gabrielov,  Degrees of real Wronski maps, 
to appear in {\it Discrete and  Computational  Geometry}.  

\bibitem[F]{F} W.~Fulton, ``Young Tableaux'', Cambridge University Press, 1997.

\bibitem[Ga]{Ga} A.~Gabrielov, Generating functions of Wronski maps,
{\it private communication}.

\bibitem[G]{G} M.~Gaudin, Diagonalization d'une class hamiltoniens de spin.
{\it Journ. de Physique} {\bf 37}, no. 10 (1976), 1087 -- 1098.       

\bibitem[GH]{GH} Ph.~Griffits and J.~Harris, 
``Principles of Algebraic Geometry'', Wiley, New York, 1978.                 

\bibitem[KZ]{KZ} V.~Knizhnik and A.~Zamolodchikov, Current algebra and 
Wess--Zumino models in two dimensions, {\it Nucl. Phys.} {\bf B247}
(1984), 83--103.

\bibitem[R]{R} E.~Rainville, Intermediate differential equations,
The Macmillan Company, 1964.

\bibitem[RV]{RV} N.~Reshetikhin and A. ~Varchenko, Quasiclassical Asymptotics
of Solutions to the KZ Equations. In: Geometry, Topology, and Physics
for Raoul Bott, International Press, 1994, 293--322.

\bibitem[SV]{SV} V.~ Schechtman and A. ~Varchenko,
Arrangements of hyperplanes and Lie algebra homology.
 {\it Invent. Math.} {\bf 106} (1991), 139 -- 194.

\bibitem[ScV]{ScV} I.~Scherbak and A.~Varchenko,  Critical points of functions,
$sl_2$ representations, and Fuchsian differential equations
with only univalued solutions, {\it preprint} (2001),  math.QA/0112269. 

\bibitem[S]{S} G.~Szego, ``Orthogonal polynomials'',  AMS, 1939.

\bibitem[V]{V} A.~Varchenko,  Special functions, KZ type equations
and Representation theory.  {\it Notes of a course given at MIT during
the spring of 2002},  math.QA/0205313.

\end{thebibliography}
\end{document}